# On Inclined Curves According to Parallel Transport Frame in $E^4$

Fatma GÖKÇELİK, İsmail GÖK, F. Nejat EKMEKCİ, and Yusuf YAYLI

**Abstract.** In this paper, we introduce an inclined curves according to parallel transport frame. Also, we define a vector field $D$ called Darboux vector field of an inclined curve in $E^4$ and we give a new characterization such as:

"$\alpha\ :\ I \subset \mathrm{R} \to E^4$ is an inclined curve $\Leftrightarrow k_1 \int k_1 ds + k_2 \int k_2 ds + k_3 \int k_3 ds = 0$"

where $k_1$, $k_2$, $k_3$ are the principal curvature functions according to parallel transport frame of the curve $\alpha$ and we give the similar characterizations such as

"$\alpha\ :\ I \subset \mathrm{R} \to E^3$ is a generalized helix $\Leftrightarrow k_1 \int k_1 ds + k_2 \int k_2 ds = 0$"

where $k_1$, $k_2$ are the principal curvature functions according to Bishop frame of the curve $\alpha$. Moreover, we illustrate some examples and draw their figures with Mathematica Programme.

## 1. Introduction

The curves are a part of our lives are the indispensable. For example, heart chest film with X-ray curve, how to act is important to us. Curves give the movements of the particle in Physics. Helical curves are very important type of curves. Because, helices are among the simplest objects in the art, molecular structures, nature, etc. For example, the path, arroused by the climbing of beans and the orbit where the progressing of the screw are a helix curves. Also, in medicine DNA molecule is formed as two intertwined helices and many proteins have helical structures, known as alpha helices. So, such curves are very important for understand to nature. Therefore, lots of author interested in the helices and they published many papers in Euclidean $3-$ and $4-$ space (See for details: [1], [3], [7], [13], [16], [17], [18], [21], [24]).

Helix curve is defined by the property that the tangent vector field makes a constant angle with a fixed direction. In 1802, M. A. Lancert first proposed a theorem and in 1845, B. de Saint Venant first proved this theorem: "A necessary and sufficient condition that a curve be a general helix is that the ratio of curvature to torsion be constant" [23].

Recently, many studies have been reported on generalized helices and inclined curves (Generalized helix is called as inclined curve in $n-$dimensional Euclidean space $E^n$, $n \geq 4$) [1], [4], [8], [12], [14], [20], [22].

The Frenet frame is constructed for the curve of $3$-time continuously differentiable non-degenerate curves. Curvature of the curve may vanish at some points on the curve, that is, second derivative of the curve may be zero. In this situation, we need an alternative frame in $E^3$. Therefore in [2], Bishop defined a new frame for a curve and called it Bishop frame which is well defined even when the curve has vanishing second derivative in $3-$dimensional Euclidean space $E^3$.

Similarly, Gökçelik et al. defined a new frame for a curve and they called it parallel transport frame in $E^4$ [9]. The parallel transport frame is an alternative frame defined a moving frame. In [9], they consider a regular curve $\alpha(s)$ parametrized by $s$ and they defined a normal vector field $V(s)$ which is perpendicular to the tangent vector field $T(s)$ of the curve $\alpha(s)$ said to be relatively parallel vector field if its derivative is tangential along the curve $\alpha(s)$. They use the tangent vector $T(s)$ and three relatively parallel vector fields to construct this alternative frame. They choose any convenient arbitrary basis $\{M_1(s), M_2(s), M_3(s)\}$ of the frame, which are perpendicular to $T(s)$ at each point. The derivatives of $\{M_1(s), M_2(s), M_3(s)\}$ only depend on $T(s)$. It is called as parallel transport frame along a curve because the normal component of the derivatives of the normal vector field is zero. The advantages of the parallel frame and the comparable parallel frame with the Frenet frame in $3-$ dimensional Euclidean space $E^3$ was given and studied by Bishop [2].

Spherical curves are characterized by the parallel transport frame. "The curve $\alpha$ is spherical curve if and only if $c_1 k_1 + c_2 k_2 + c_2 k_3 + 1 = 0$, where $c_1$, $c_2$, $c_3$ are constant and $k_1$, $k_2$, $k_3$ are the principal curvatures according to parallel transport frame of the curve $\alpha$ " [9].

In this article, we study inclined curves according to parallel transport frame in terms of the harmonic curvature functions and give a characterization such as: "The curve $\alpha : I \subset \mathrm{R} \to E^4$ is a generalized helix if and only if $k_1 \int k_1 ds + k_2 \int k_2 ds + k_3 \int k_3 ds = 0$ and $k_1$, $k_2$, $k_3$ show the principal curvature functions according to parallel transport frame of the curve $\alpha$ ".

## 2. Preliminaries

Let $\alpha : I \subset \mathrm{R} \to E^4$ be an arbitrary curve in the Euclidean $4-$ space $E^4$. Recall that the curve $\alpha$ is said to be of unit speed (or parametrized by arc-length function $s$ ) if $\langle \alpha'(s), \alpha'(s) \rangle = 1$, where $\langle , \rangle$ is the standard inner product of $E^4$ given by

$$\langle X, Y \rangle = x_1 y_1 + x_2 y_2 + x_3 y_3 + x_4 y_4$$

for each $X = (x_1, x_2, x_3, x_4)$, $Y = (y_1, y_2, y_3, y_4) \in E^4$. In particular, the norm of a vector $X \in E^4$ is given by $\|X\|^2 = \langle X, X \rangle$. Let $\{T, N, B_1, B_2\}$ be the Frenet frame along the unit speed curve $\alpha$. Then $T, N, B_1$ and $B_2$ are the tangent, the principal normal, first and second binormal vectors of the curve $\alpha$, respectively. If
$\alpha$ is a space curve, then this set of orthogonal unit vectors, known as the Frenet-Serret frame, has the following properties

$$\begin{aligned} T'(s) &= \bar{k}_1 N(s) \\ N'(s) &= -\bar{k}_1 T(s) + \bar{k}_2 B(s) \\ B_1'(s) &= -\bar{k}_2 N(s) + \bar{k}_3 B_2(s) \\ B_2'(s) &= -\bar{k}_3 B_1(s) \end{aligned} \quad (1.1)$$

where $\bar{k}_1$, $\bar{k}_2$ and $\bar{k}_3$ denote principal curvature functions according to Frenet frame of the curve $\alpha$ [10].

The parallel transport frame is an alternative frame defined a moving frame. Curvature of the curve may vanish at some points on the curve, that is, the $i-th$ $(1<i<4)$ derivative of the curve may be zero. We can parellel transport an orthonormal frame along a curve simply by parallel transporting each component of the frame. The derivatives of $\{M_1(s), M_2(s), M_3(s)\}$ only depend on $T(s)$.

Here the set $\{T(s), M_1(s), M_2(s), M_2(s)\}$ is called as parallel transport frame and
$k_1(s) = \langle T'(s), M_1(s) \rangle$, $k_2(s) = \langle T'(s), M_2(s) \rangle$, $k_3 = \langle T'(s), B_2(s) \rangle$ called as parallel transport curvatures of the curve $\alpha$.

**Theorem** *Let* $\{T, N, B_1, B_2\}$ *be a Frenet frame along a unit speed curve* $\alpha : I \subset IR \to E^4$ *and* $\{T, M_1, M_2, M_3\}$ *denotes the parallel transport frame of the curve* $\alpha$. *The relation may be expressed as*

$$T = T,$$
$$N = \cos\theta(s)\cos\psi(s)M_1 + (-\cos\phi(s)\sin\psi(s) + \sin\phi(s)\sin\theta(s)\cos\psi(s))M_2$$
$$+(\sin\phi(s)\sin\psi(s) + \cos\phi(s)\sin\theta(s)\cos\psi(s))M_3,$$

$$B_1 = \cos\theta(s)\sin\psi(s)M_1 + (\cos\phi(s)\cos\psi(s) + \sin\phi(s)\sin\theta(s)\sin\psi(s))M_2$$
$$+(-\sin\phi(s)\cos\psi(s) + \cos\phi(s)\sin\theta(s)\sin\psi(s))M_3,$$

$$B_2 = -\sin\theta(s)M_1 + \sin\phi(s)\cos\theta(s)M_2 + \cos\phi(s)\cos\theta(s)M_3$$

*where* $\theta$, $\phi$, $\psi$ *are angles between vectors of parallel transport frame which are shown Figure 1* [9].

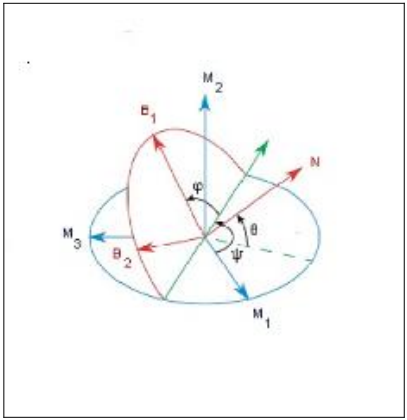

**Figure 1** : The relation of the Frenet and the Parallel Transport frame by means of the rotation matrix.

**Theorem 2.** *The alternative parallel transport frame equations are given by*

$$\begin{bmatrix} T' \\ M_1' \\ M_2' \\ M_3' \end{bmatrix} = \begin{bmatrix} 0 & k_1 & k_2 & k_3 \\ -k_1 & 0 & 0 & 0 \\ -k_2 & 0 & 0 & 0 \\ -k_3 & 0 & 0 & 0 \end{bmatrix} \begin{bmatrix} T \\ M_1 \\ M_2 \\ M_3 \end{bmatrix} \quad (2.1)$$

*where* $k_1, k_2, k_3$ *are principal curvature functions acoording to parallel transport frame of the curve* $\alpha$ *and their expression as follows:*

$$\begin{aligned} k_1 &= \bar{k}_1 \cos\theta \cos\psi, \\ k_2 &= \bar{k}_1(-\cos\phi\sin\psi + \sin\phi\sin\theta\cos\psi), \\ k_3 &= \bar{k}_1(\sin\phi\sin\psi + \cos\phi\sin\theta\cos\psi) \end{aligned} \quad (2.2)$$

$$\begin{aligned} \bar{k}_1(s) &= \sqrt{k_1^2 + k_2^2 + k_3^2}, \\ \bar{k}_2(s) &= -\psi' + \phi'\sin\theta, \\ \bar{k}_3(s) &= \frac{\theta'}{\sin\psi}, \end{aligned} \quad (2.3)$$

$$\phi'\cos\theta + \theta'\cot\psi = 0.$$

*where* $\bar{k}_1, \bar{k}_2, \bar{k}_3$ *are the principal curvature functions according to Frenet frame of the curve* $\alpha$ [9].

If we take $\bar{k}_3(s) = 0$ then we get Bishop frame [2] in $3-$ dimensional Euclidean space $E^3$ defined as:

$$\begin{bmatrix} T' \\ M_1' \\ M_2' \end{bmatrix} = \begin{bmatrix} 0 & k_1 & k_2 \\ -k_1 & 0 & 0 \\ -k_2 & 0 & 0 \end{bmatrix} \begin{bmatrix} T \\ M_1 \\ M_2 \end{bmatrix}.$$

# 3. On Inclined Curves According to Parallel Transport Frame in $E^4$

In this section, we give some characterizations for an inclined curve parallel transport frame in $E^4$ by using the harmonic curvature functions of the curve.

**Definition 1.** *Let* $\alpha : I \subset \mathrm{R} \to E^4$ *be a unit speed curve in* $4-$ *dimensional Euclidean space* $E^4$ *with the arc-length parameter* $s$ *and* $X$ *be a unit and fixed vector of* $E^4$. *For all* $s \in I,$ *if*

$$\langle \alpha'(s), X \rangle = \cos\varphi = \text{const.}, \varphi \neq \frac{\pi}{2}$$

then the curve $\alpha$ is called an inclined curve in $E^4$ where $\alpha'(s)$ is the unit tangent vector field to the curve $\alpha$ at its point $\alpha(s)$ and $\varphi$ is a constant angle between the vectors $\alpha'$ and $X$.

**Definition 2.** Let $\alpha : I \subset R \to E^4$ be a regular curve with arc-length parameter $s$ and $\{k_1, k_2, k_3\}$ are nonzero curvatures according to parallel transport frame. In this case the harmonic curvature functions of the curve $\alpha$ are given by

$$H_i : I \to R, (i = 1, 2, 3)$$
$$H_1(s) = -[AH_2(s) + BH_3(s)],$$

$$H_2(s) = -\frac{B'}{A'} \frac{1}{(\frac{B'}{A'})'} [k_2 + \frac{B'}{A'} k_3 + (\frac{k_1^2 + k_2^2 + k_3^2}{k_1 A})'] + \frac{k_1^2 + k_2^2 + k_3^2}{k_1 A},$$

$$H_3(s) = \frac{1}{(\frac{B'}{A'})'} [k_2 + \frac{B'}{A'} k_3 + (\frac{k_1^2 + k_2^2 + k_3^2}{k_1 A})']$$

where $A = \frac{k_2}{k_1}$, $B = \frac{k_3}{k_1}$.

**Theorem 3.** Let $\alpha : I \subset R \to E^4$ be a regular curve with arc-length parameter $s$ and $\{k_1, k_2, k_3\}$ are nonzero curvatures according to parallel transport frame and $\{T(s), M_1(s), M_2(s), M_3(s)\}$ denotes the parallel transport frame of the curve $\alpha$. If the curve $\alpha$ is an inclined curve in $E^4$ then

$$\langle M_1, X \rangle = H_1 \langle T, X \rangle,$$
$$\langle M_2, X \rangle = H_2 \langle T, X \rangle, \quad (3.1)$$
$$\langle M_3, X \rangle = H_3 \langle T, X \rangle$$

where $H_i(s)$ $(i = 1, 2, 3)$ are the harmonic curvature functions of the curve $\alpha$.

**Proof.** Let $\alpha$ be an inclined curve in $E^4$. Then
$$\langle T, X \rangle = \cos\varphi = \text{constant}.$$

Differentiating the above equation, with respect to $s$, we obtain
$$\langle T', X \rangle = 0,$$
$$\langle k_1 M_1 + k_2 M_2 + k_3 M_3, X \rangle = 0,$$
$$k_1 \langle M_1, X \rangle + k_2 \langle M_2, X \rangle + k_3 \langle M_3, X \rangle = 0,$$

$$\langle M_1, X \rangle = -[A \langle M_2, X \rangle + B \langle M_3, X \rangle].$$

Again differentiating of the last equation

$$\langle M_2, X \rangle = -\frac{B'}{A'} \langle M_3, X \rangle + (\frac{k_1^2 + k_2^2 + k_3^2}{k_1 A}) \langle T, X \rangle,$$

$$\langle M_3, X \rangle = \frac{1}{(\frac{B'}{A'})'}[k_2 + \frac{B'}{A'}k_3 + (\frac{k_1^2 + k_2^2 + k_3^2}{k_1 A})']\langle T, X \rangle.$$

*Then by using the* Definition (2) *it is easy to obtain*

$$\langle M_1, X \rangle = H_1 \langle T, X \rangle,$$
$$\langle M_2, X \rangle = H_2 \langle T, X \rangle,$$
$$\langle M_3, X \rangle = H_3 \langle T, X \rangle.$$

Theorem 4. *Let* $\alpha : I \subset \mathrm{R} \to E^4$ *be a regular curve with arc-length parameter* $s$ *and* $\{T(s), M_1(s), M_2(s), M_3(s)\}$ *denotes the parallel transport frame of the curve* $\alpha$. *If the curve* $\alpha$ *is an inclined curve in* $E^4$, *the axis of the curve* $\alpha$ *is* $X$ *given by*

$$X = (T(s) + H_1 M_1(s) + H_2 M_2(s) + H_3 M_3(s))\langle T(s), X \rangle$$

*or*

$$X = (T(s) + H_1 M_1(s) + H_2 M_2(s) + H_3 M_3(s))\cos\varphi.$$

*where* $H_i(s)$, $(i = 1, 2, 3)$ *are the harmonic curvature functions of the curve* $\alpha$.

Proof. If the axis of the inclined curve $\alpha$ is $X$, then we get

$$X = \lambda_1 T + \lambda_2 M_1 + \lambda_3 M_2 + \lambda_4 M_3.$$

Since $\alpha$ is an inclined curve $\lambda_1 = \cos\varphi$ and then by using Theorem (3)

$$\lambda_2 = H_1 \langle T, X \rangle,$$
$$\lambda_3 = H_2 \langle T, X \rangle,$$
$$\lambda_4 = H_3 \langle T, X \rangle.$$

Thus it is easy to obtain

$$X = (T(s) + H_1 M_1(s) + H_2 M_2(s) + H_3 M_3(s))\cos\varphi.$$

Definition 3. *Let* $\alpha : I \subset \mathrm{R} \to E^4$ *be a regular curve with arc-length parameter* $s$ *and* $\{T(s), M_1(s), M_2(s), M_3(s)\}$ *denotes the parallel transport frame of the curve* $\alpha$ *and* $H_i$ $(i = 1, 2, 3)$ *denote the harmonic curvature functions at the point* $\alpha(s)$.

$$D = T(s) + H_1 M_1(s) + H_2 M_2(s) + H_3 M_3(s)$$

*is called a Darboux vector field of an inclined curve* $\alpha$ *in* $E^4$.

Theorem 5. *Let* $\alpha : I \subset \mathrm{R} \to E^4$ *be a regular curve with arc-length parameter* $s$ *and* $\{T(s), M_1(s), M_2(s), M_3(s)\}$ *denotes the parallel transport frame of the curve* $\alpha$ *and* $H_i$ $(i = 1, 2, 3)$ *denote the harmonic curvature functions at the point* $\alpha(s)$. *Then the curve* $\alpha$ *is an inclined curve in* $E^4$ *if and only if* $D$ *is a constant vector field.*

Proof. Let $\alpha$ be an inclined curve in $E^4$ and $X$ is the axis of the curve $\alpha$. From Theorem (4), we have

$$X = (T(s) + H_1 M_1(s) + H_2 M_2(s) + H_3 M_3(s))\cos\varphi$$
$$= D\cos\varphi$$

where $\cos\varphi$ and $X$ are constant. Hence $D$ be a constant vector field.

Conversely, suppose that $D$ be a constant vector field, then we have $\|D\|^2 = \langle D, D\rangle$ is constant. By using Theorem (4) we get

$$\|X\|^2 = \|D\cos\varphi\|^2$$
$$= \langle D\cos\varphi, D\cos\varphi\rangle$$
$$= \cos^2\varphi\langle D, D\rangle.$$

Since $X$ is a unit vector and $\|D\|$ is constant, we have $\cos\varphi = \frac{1}{\|D\|} = \langle T(s), X\rangle$ is a constant. So, $\alpha$ is an inclined curve. This completes the proof.

**Theorem 6.** Let $\alpha : I \subset \mathbb{R} \to E^4$ be a regular curve with arc-length parameter $s$ and $\{T(s), M_1(s), M_2(s), M_3(s)\}$ denotes the parallel transport frame of the curve $\alpha$ and $H_i$ $(i=1,2,3)$ denote the harmonic curvature functions at the point $\alpha(s)$. If the curve $\alpha$ is an inclined curve in $E^4$ then $\sum_{i=1}^{3} H_i^2(s)$ is a constant.

**Proof.** Let $\alpha$ be an inclined curve. Since, the axis of the curve $\alpha$ is
$$X = (T(s) + H_1 M_1(s) + H_2 M_2(s) + H_3 M_3(s))\cos\varphi = D\cos\varphi \quad \text{unit vector field,} \quad \|X\|^2 = 1.$$
Hence, from the Theorem (4)

$$\|X\|^2 = \langle X, X\rangle$$
$$= \cos^2\varphi + \sum_{i=1}^{3} H_i^2(s)\cos^2\varphi$$

we obtain

$$\|X\|^2 = 1 = (1 + \sum_{i=1}^{3} H_i^2(s))\cos^2\varphi$$

and then

$$\sum_{i=1}^{3} H_i^2(s) = \tan^2\varphi = \text{constant.}$$

**Theorem 7.** Let $\alpha : I \subset \mathbb{R} \to E^4$ be a curve with arc-length parameter $s$ and $\{k_1, k_2, k_3\}$ are nonzero curvatures according to parallel transport frame. Such that $\{T(s), M_1(s), M_2(s), M_3(s)\}$ denotes the parallel transport frame of the curve $\alpha$, $\{H_1, H_2, H_3\}$ denotes the harmonic curvature functions of the curve $\alpha$. Then $\alpha$ is an inclined curve if and only if

$$H_1' = -k_1 \ , \ H_2' = -k_2 \ \text{and} \ H_3' = -k_3. \tag{3.2}$$

**Proof** If we differentiate $D$ along the curve $\alpha$, we get

$$D' = T'(s) + H_1'M_1(s) + H_1M_1'(s) + H_2'M_2(s) + H_2M_2'(s)$$
$$+ H_3'M_3(s) + H_3M_3'(s).$$

Using the alternative parallel frame equations we have

$$D' = (k_1 + H_1')M_1 + (k_2 + H_2')M_2 + (k_3 + H_3')M_3.$$

From Theorem (5) since $\alpha$ is an inclined curve, $D$ is a constant vector field. So, the Eq. (3.2) holds.

Conversely, if the Eq. (3.2) holds, we can easily see that $D' = 0$ or $D$ is a constant vector field, and then from Theorem (5), we have that $\alpha$ is an inclined curve in $E^4$. This completes the proof.

**Corollary 1.** *Let* $\alpha : I \subset R \to E^4$ *be a curve with arc-length parameter* $s$ *and* $\{k_1, k_2, k_3\}$ *are nonzero curvatures according to parallel transport frame. Such that* $\{T(s), M_1(s), M_2(s), M_3(s)\}$ *denotes the parallel transport frame of the curve* $\alpha$, $\{H_1, H_2, H_3\}$ *denotes the harmonic curvature functions of the curve* $\alpha$. *Then* $\alpha$ *is an inclined curve if and only if*

$$k_1(s)\int k_1(s)ds + k_2(s)\int k_2(s)ds + k_3(s)\int k_3(s)ds = 0. \qquad (3.3)$$

**Proof.** Let $\alpha$ be an inclined curve in $E^4$, then from Theorem (7) we have $H_1' = -k_1$, $H_2' = -k_2$ and $H_3' = -k_3$. By using Theorem (7)

$$k_1(s)\int k_1(s)ds + k_2(s)\int k_2(s)ds + k_3(s)\int k_3(s)ds = 0.$$

Conversely, we suppose that the Eq. (3.3) holds. From Theorem (7) and Definition (2), it is obvious that $\alpha$ is an inclined curve in $E^4$.

**Definition 4.** *The matrix* $F_4(s)$ *with respect to the basis* $T(s), M_1(s), M_2(s), M_3(s)$ *is given by*

$$F_4(s) = \begin{bmatrix} 0 & k_1 & k_2 & k_3 \\ -k_1 & 0 & 0 & 0 \\ -k_2 & 0 & 0 & 0 \\ -k_3 & 0 & 0 & 0 \end{bmatrix} \in R_4^4.$$

**Corollary 2.** *Let* $\alpha$ *be a unit speed non-degenerate curve in* $E^4$, $\{k_1, k_2, k_3\}$ *denotes the curvature functions according to parallel transport frame of the curve* $\alpha$ *and the matrix* $F_4(s)$. *The curve* $\alpha$ *is an inclined curve if and only if the vector* $D = [1, H_1, H_2, H_3] \in R^4$ *satisfies equations*

$$\frac{d}{ds}\begin{bmatrix} 1 \\ H_1 \\ H_2 \\ H_3 \end{bmatrix} = F_4(s)\begin{bmatrix} 1 \\ H_1 \\ H_2 \\ H_3 \end{bmatrix}.$$

**Proof.** Direct substitution shows that

$$F_4(s)D^T = \begin{bmatrix} 0 & k_1 & k_2 & k_3 \\ -k_1 & 0 & 0 & 0 \\ -k_2 & 0 & 0 & 0 \\ -k_3 & 0 & 0 & 0 \end{bmatrix}\begin{bmatrix} 1 \\ H_1 \\ H_2 \\ H_3 \end{bmatrix}$$

$$= [0, H'_1, H'_2, H'_3]^T$$

Since $\alpha$ is an inclined curve

$$\frac{d}{ds}[1, H_1, H_2, H_3]^T = F_4(s)[1, H_1, H_2, H_3]^T.$$

**Theorem 8.** *Let $\alpha$ be a unit speed non-degenerate curve in $E^4$. $\{k_1, k_2, k_3\}$ denotes the curvature functions and $\{T(s), M_1(s), M_2(s), M_3(s)\}$ denotes the parallel transport frame of the curve $\alpha$. $\alpha_{M_1}, \alpha_{M_2}, \alpha_{M_3}$ are the spherical images of the $M_1$, $M_2$, $M_3$, respectively. $\alpha$ is an inclined curve iff the spherical image curves $\alpha_{M_1}, \alpha_{M_2}, \alpha_{M_3}$ are inclined curve. Also, axis of these curves are the same axis of the curve $\alpha$.*

**Proof.** Let $\alpha$ be a unit speed non-degenerate curve in $E^4$. Its spherical image is the parametrized curve with arc-length parametes $s_{M_1}$, then we have

$$\alpha_{M_1}(s_{M_1}) = M_1(s).$$

Differentiating the above equation with respect to $s$ we obtain

$$\frac{d\alpha_{M_1}}{ds_{M_1}}\frac{ds_{M_1}}{ds} = \frac{dM_1}{ds}.$$

If we consider then using the Eq. (2.1) we get

$$T_{M_1} = (-k_1 T)\frac{ds}{ds_{M_1}}$$

and taking the norm of the tangent vector field $T_{M_1}$ of the curve $\alpha_{M_1}$ we get

$$\frac{ds_{M_1}}{ds} = k_1.$$

Thus $T_{M_1} = T$. Consequently, if $\langle T, X \rangle$ is a constant, then $\langle T_{M_1}, X \rangle$, $\langle T_{M_2}, X \rangle$ and $\langle T_{M_3}, X \rangle$ must be constant. Hence, $\alpha$ is an inclined curve in $E^4$ then $\alpha_{M_1}, \alpha_{M_2}, \alpha_{M_3}$ are the inclined curve.

In this section, we will give some theorems without their proofs about the generalized helices in $E^3$ according to Bishop frame. Because their proofs can be shown using the similar method of the above proofs in inclined curves in $E^4$. Thus, we manage to give some examples of the generalized helices according to Bishop frame and draw their figures by using Mathematica Programme.

**Definition 5.** *Let $\alpha : I \to E^3$ be a unit speed curve in $E^3$ with the arc-length parameter $s$ and $X$ be a unit constant vector of $E^3$. For all $s \in I$, if*

$$\langle \alpha'(s), X \rangle = \cos\theta = \text{const.}, \theta \neq \frac{\pi}{2}$$

*then the curve $\alpha$ is called a generalized helix in $E^3$ where $\alpha'(s)$ is the unit tangent vector field to the curve $\alpha$ at its point $\alpha(s)$ and $\theta$ is a constant angle between the vectors $\alpha'$ and $X$.*

**Definition 6.** *Let $\alpha : I \subset \mathrm{R} \to E^3$ be a regular curve with arc-length parameter $s$ and $\{k_1, k_2\}$ are nonzero curvatures according to Bishop frame . In that case, the harmonic curvature functions of the curve $\alpha$ are given by*

$$H_i : I \to \mathrm{R}, (i = 1, 2)$$

$$H_1(s) = \frac{k_2(1 + f^2)}{f'},$$

$$H_2(s) = -\frac{k_1(1 + f^2)}{f'}$$

*where $f = \frac{k_1}{k_2}$.*

**Theorem 9.** *Let $\alpha : I \subset \mathrm{R} \to E^3$ be a regular curve with arc-length parameter $s$ and $\{k_1, k_2\}$ are nonzero curvatures according to Bishop frame and $\{T(s), M_1(s), M_2(s)\}$ denotes the Bishop frame of the curve $\alpha$. If the curve $\alpha$ is a generalized helix in $E^3$ then*

$$\langle M_1, X \rangle = H_1 \langle T, X \rangle,$$

$$\langle M_2, X \rangle = H_2 \langle T, X \rangle$$

*where $H_i(s)$ $(i = 1, 2)$ are harmonic curvature functions of the curve $\alpha$.*

**Theorem 10.** *Let $\alpha : I \subset \mathrm{R} \to E^3$ be a regular curve with arc-length parameter $s$ and $\{T(s), M_1(s), M_2(s)\}$ denotes the Bishop frame of the curve $\alpha$. If the curve $\alpha$ is generalized helix in $E^3$, the axis of $\alpha$ is $X$*

$$X = (T(s) + H_1 M_1(s) + H_2 M_2(s)) \langle T(s), X \rangle$$

*or*

$$X = (T(s) + H_1 M_1(s) + H_2 M_2(s)) \cos\theta.$$

*where $H_i(s)$ $(i = 1, 2)$ are harmonic curvature functions of the curve $\alpha$.*

**Definition 7.** *Let $\alpha : I \subset \mathrm{R} \to E^3$ be a regular curve with arc-length parameter $s$ and $\{T(s), M_1(s), M_2(s)\}$ denotes the Bishop frame of the curve $\alpha$ and $H_i$ $(i = 1, 2)$ denote the harmonic curvature functions at the point $\alpha(s)$.*

$$D = T(s) + H_1 M_1(s) + H_2 M_2(s)$$

is called a Darboux vector field generalized helix $\alpha$ in $E^3$.

**Theorem 11.** *Let* $\alpha : I \subset R \to E^3$ *be a regular curve with arc-length parameter* $s$ *and* $\{T(s), M_1(s), M_2(s)\}$ *denotes the Bishop frame of the curve* $\alpha$ *and* $H_i$ $(i=1,2)$ *denote the harmonic curvature functions at the point* $\alpha(s)$. *Then the curve* $\alpha$ *is generalized helix in* $E^3$ *if and only if* $D$ *is a constant vector field.*

**Theorem 12.** *Let* $\alpha : I \subset R \to E^3$ *be a regular curve with arc-length parameter* $s$ *and* $\{T(s), M_1(s), M_2(s)\}$ *denotes the Bishop frame of the curve* $\alpha$ *and* $H_i$ $(i=1,2)$ *denote the harmonic curvature functions at the point* $\alpha(s)$. *Then the curve* $\alpha$ *is a generalized helix in* $E^3$ *if and only if* $H_1^2(s) + H_2^2(s)$ *is a constant.*

**Theorem 13.** *Let* $\alpha : I \subset R \to E^3$ *be a curve with arc-length parameter* $s$ *and* $\{k_1, k_2\}$ *are nonzero curvatures according to Bishop frame. Such that* $\{T(s), M_1(s), M_2(s)\}$ *denotes the Bishop frame of the curve* $\alpha$, $\{H_1, H_2\}$ *denotes the harmonic curvature functions of the curve* $\alpha$. *Then* $\alpha$ *is a generalized helix if and only if*

$$H_1' = -k_1 \text{ and } H_2' = -k_2.$$

**Corollary 3.** *Let* $\alpha : I \subset R \to E^3$ *be a curve with arc-length parameter* $s$ *and* $\{k_1, k_2\}$ *are nonzero curvatures according to Bishop frame. Such that* $\{T(s), M_1(s), M_2(s)\}$ *denotes the Bishop frame and* $\{H_1, H_2\}$ *denotes the harmonic curvature functions of the curve* $\alpha$. *Then* $\alpha$ *is a generalized helix if and only if*

$$k_1(s)\int k_1(s)ds + k_2(s)\int k_2(s)ds = 0.$$

**Lemma 1.** *The matrix* $F_3(s)$ *with respect to the basis* $T(s), M_1(s), M_2(s)$ *is given by*

$$F_3(s) = \begin{bmatrix} 0 & k_1 & k_2 \\ -k_1 & 0 & 0 \\ -k_2 & 0 & 0 \end{bmatrix} \in R_3^3.$$

**Corollary 4.** *Let* $\alpha$ *be a unit speed non-degenerate curve in* $E^3$, $\{k_1, k_2\}$ *denotes the curvature functions according to Bishop frame of the curve* $\alpha$ *and the matrix* $F_3(s)$. *The curve* $\alpha$ *is a generalized helix if and only if the vector* $D = [1, H_1, H_2] \in R^3$ *satisfies the Frenet equations:*

$$\frac{d}{ds}\begin{bmatrix} 1 \\ H_1 \\ H_2 \end{bmatrix} = F_3(s)\begin{bmatrix} 1 \\ H_1 \\ H_2 \end{bmatrix}.$$

## 4. Examples

In this section, we illustrate two examples of generalized helices according to Bishop frame and we draw their figures by using the Mathematica Programme.

**Example 1.** *Let us consider the following curve* $\beta(s) = (\beta_1, \beta_2, \beta_3)$ *of* $E^3$

$$\begin{cases} \beta_1 = \cos s \cos \sqrt{17}s + \frac{1}{\sqrt{17}} \sin s \sin \sqrt{17}s \\ \beta_2 = -\cos s \sin \sqrt{17}s + \frac{1}{\sqrt{17}} \sin s \cos \sqrt{17}s \\ \beta_3 = \frac{4}{\sqrt{17}} \sin s \end{cases}$$

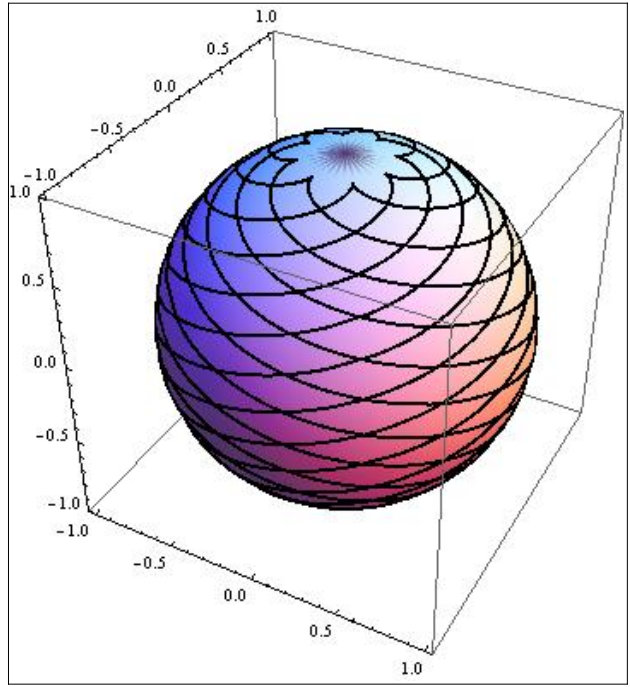

Figure 2 : The generalized helix $\beta = \beta(s)$.

It is rendered in Figure 2.
And, this curve's curvature functions are calculated with the help of Mathematica Programme.

$$\kappa = \frac{1}{\cos s}$$

$$\tau = -\frac{1}{4\cos s}.$$

The Frenet-Serret frame of the curve $\beta = \beta(s)$ may be written by the aid Mathematica Programme as follows

$$T = \left(-\frac{4}{\sqrt{17}}\sin\sqrt{17}s, -\frac{4}{\sqrt{17}}\cos\sqrt{17}s, \frac{1}{\sqrt{17}}\right),$$

$$N = (-\cos\sqrt{17}s, \sin\sqrt{17}s, 0),$$

$$B = (-\frac{1}{\sqrt{17}}\sin\sqrt{17}s, -\frac{1}{\sqrt{17}}\cos\sqrt{17}s, -\frac{4}{\sqrt{17}}).$$

To create a Bishop frame to find the angle of rotation $\beta = \beta(s)$ has the form

$$\theta(s) = -\int \tau ds = \frac{1}{4}\ln(\frac{1+\sin s}{\cos s}).$$

The transformation matrix for the curve $\beta$

$$\begin{bmatrix} T \\ N \\ B \end{bmatrix} = \begin{bmatrix} 1 & 0 & 0 \\ 0 & \cos(\frac{1}{4}\ln(\frac{1+\sin s}{\cos s})) & -\sin(\frac{1}{4}\ln(\frac{1+\sin s}{\cos s})) \\ 0 & \sin(\frac{1}{4}\ln(\frac{1+\sin s}{\cos s})) & \cos(\frac{1}{4}\ln(\frac{1+\sin s}{\cos s})) \end{bmatrix} \begin{bmatrix} T \\ M_1 \\ M_2 \end{bmatrix}$$

$T, N_1, N_2$ can easily be found

$$T(s) = T(s),$$

$$M_1(s) = \cos(\frac{1}{4}\ln(\frac{1+\sin s}{\cos s}))N(s) + \sin(\frac{1}{4}\ln(\frac{1+\sin s}{\cos s}))B(s),$$

$$M_2(s) = -\sin(\frac{1}{4}\ln(\frac{1+\sin s}{\cos s}))N(s) + \cos(\frac{1}{4}\ln(\frac{1+\sin s}{\cos s}))B(s).$$

Even, first curvature function according to Bishop frame of the curve $\beta$ is calculated

$$k_1(s) = \langle T'(s), M_1(s) \rangle = \frac{1}{\cos s}\cos(\frac{1}{4}\ln(\frac{1+\sin s}{\cos s}))$$

and

$$\kappa(s) = \sqrt{k_1^2(s) + k_2^2(s)}$$

then

$$k_2(s) = \frac{1}{\cos s}\sin(\frac{1}{4}\ln(\frac{1+\sin s}{\cos s})).$$

We can see that it clearly

$$k_1(s)\int k_1(s)ds + k_2(s)\int k_2(s)ds = 0.$$

Then, the curve $\beta = \beta(s)$ is a generalized helix.

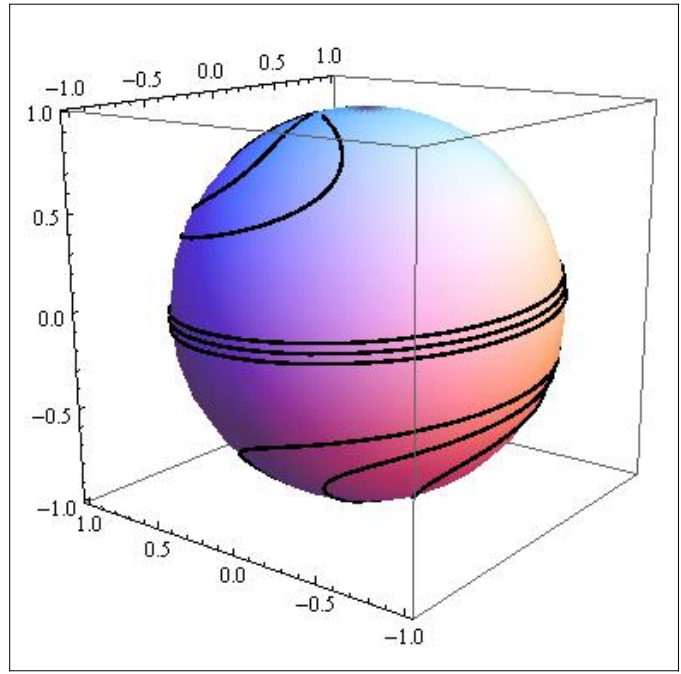

**Figure 3** : $M_1$ Bishop Spherical Images of $\beta = \beta(s)$.

**Example 2.** *Euler Spirals, were discovered indepently by three researchers. In 1694, Bernoulli wrote the equations for the Euler spiral for the first time, but did not draw the spirals or compute them numerically. In 1744, Euler rediscovered the curve's equations, described their properties, and derived a series expansion to the curve's integrals. Later, in 1781, he also computed the spiral's end points. The curves were re-discovered in 1890 for the third time by Talbot, who used them to design railway tracks [7]. An Euler spiral is a curve whose curvature changes linearly with its curve length (the curvature of a circular curve is equal to the reciprocal of the radius). Euler spirals are also commonly referred to as spiros, clothoids or Cornu spirals. Moreover, Euler spirals have applications to diffraction computations. They are also widely used as transition curve in railroad/highway engineering for connecting and transiting the geometry between a tangent and a circular curve. Let us consider the Euler Spiral* $\gamma(s) = (\gamma_1, \gamma_2, \gamma_3)$ *of* $E^3$

$$\begin{cases} \gamma_1(s) = \tfrac{3}{5}\int \sin(s^2 + 1)ds \\ \gamma_2(s) = \tfrac{3}{5}\int \cos(s^2 + 1)ds \\ \quad\gamma_3(s) = \tfrac{4}{5}s \end{cases}$$

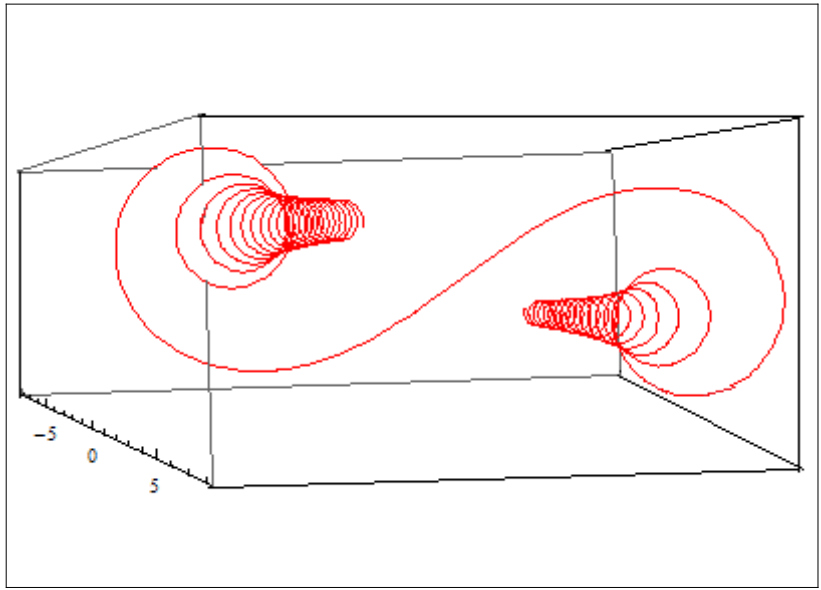

Figure 4 : The generalized helix $\gamma = \gamma(s)$.

And, we are calculated this curve's curvature functions with the help of Mathematica Programme

$$\kappa = \frac{6s}{5}$$

$$\tau = -\frac{8s}{5}.$$

The Frenet-Serret frame of the curve $\gamma = \gamma(s)$ may be written by the aid Mathematica Programme as follows

$$T(s) = \left(\frac{3}{5}\sin(s^2+1), \frac{3}{5}\cos(s^2+1), \frac{4}{5}\right),$$

$$N(s) = (\cos(s^2+1), -\sin(s^2+1), 0),$$

$$B(s) = (\frac{4}{5}\sin(s^2+1), \frac{4}{5}\cos(s^2+1), -\frac{3}{5}).$$

To create a Bishop frame to find the angle of rotation. $\gamma = \gamma(s)$ has the form

$$\theta(s) = -\int \tau ds = \frac{4s^2}{5}.$$

The transformation matrix for the curve

$$\begin{bmatrix} T \\ N \\ B \end{bmatrix} = \begin{bmatrix} 1 & 0 & 0 \\ 0 & \cos(\frac{4s^2}{5}) & -\sin(\frac{4s^2}{5}) \\ 0 & \sin(\frac{4s^2}{5}) & \cos(\frac{4s^2}{5}) \end{bmatrix} \begin{bmatrix} T \\ M_1 \\ M_2 \end{bmatrix}$$

$T, M_1, M_2$ can easily be found

$$T(s) = T(s),$$

$$M_1(s) = \cos(\frac{4s^2}{5})N(s) + \sin(\frac{4s^2}{5})B(s),$$

$$M_2(s) = -\sin(\frac{4s^2}{5})N(s) + \cos(\frac{4s^2}{5})B(s).$$

Even, first curvature function according to Bishop frame of the curve $\gamma$ is calculated

$$k_1(s) = \langle T'(s), N_1(s) \rangle = \frac{6s}{5}\cos(\frac{4s^2}{5})$$

and

$$\kappa(s) = \sqrt{k_1^2(s) + k_2^2(s)}$$

then

$$k_2(s) = \frac{6s}{5}\sin(\frac{4s^2}{5}).$$

We can see that it clearly

$$k_1(s)\int k_1(s)ds + k_2(s)\int k_2(s)ds = 0.$$

Then, the Euler Spiral $\gamma = \gamma(s)$ is generalized helix according to Bishop frame.

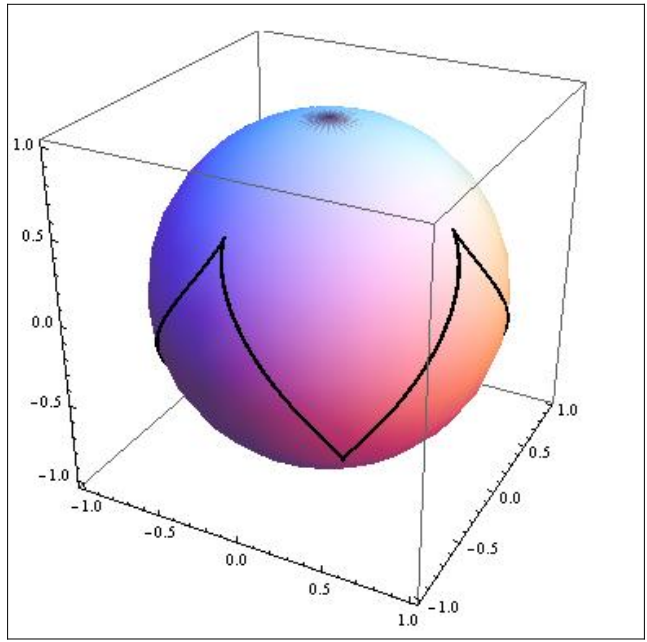

Figure 5 : $M_1$ Bishop Spherical Images of $\gamma = \gamma(s)$.

Department of Mathematics, Faculty of Science, University of Ankara Tandogan, Ankara, TURKEY
fgokcelik@ankara.edu.tr

Department of Mathematics, Faculty of Science, University of Ankara Tandogan, Ankara, TURKEY
igok@science.ankara.edu.tr

Department of Mathematics, Faculty of Science, University of Ankara Tandogan, Ankara, TURKEY
ekmekci@science.ankara.edu.tr

Department of Mathematics, Faculty of Science, University of Ankara Tandogan, Ankara, TURKEY
yyayli@science.ankara.edu.tr